\documentclass{amsart}[11]
\usepackage{amssymb}
\usepackage{graphicx}
 
\usepackage{a4} 

\usepackage{color}

\def\blue{\color{blue}}

\begin{document}
\newtheorem{theorem}{Theorem}[section]
\newtheorem{lemma}[theorem]{Lemma}
\newtheorem{proposition}[theorem]{Proposition}
\newtheorem{remark}[theorem]{Remark}
\newtheorem{definition}[theorem]{Definition}
\newtheorem{corollary}[theorem]{Corollary}
\newtheorem{example}[theorem]{Example}
\newtheorem{conjecture}[theorem]{Conjecture}
\newtheorem{problem}[theorem]{Problem}

\def\sym{\text{Sym}}
\def\qedbox{\hbox{$\rlap{$\sqcap$}\sqcup$}}

\def\CP{\mathbb{CP}}
\def\A{\mathcal{A}}
\def\w{\omega}

\makeatletter
ÊÊ\renewcommand{\theequation}{\thesection.\alph{equation}}
ÊÊ\@addtoreset{equation}{section}
\makeatother

\title[Nets in $\mathbb{CP}^2$]
{equivalence classes of latin squares and nets in $\mathbb{CP}^2$}

\author{C. Dunn}
\address{CD:  Department of Mathematics, California State University, San Bernardino, 92407.}
\email{cmdunn@csusb.edu}

\author{M. Miller}
\address{MM:  Department of Mathematics, University of Oregon, Eugene, OR, 97403.}
\email{mmiller7@uoregon.edu}

\author{M. Wakefield}
\address{MW:  Department of Mathematics, Hokkaido University, Sapporo, 060-0810, Japan.}
\curraddr{}
\email{wakefield@math.sci.hokudai.ac.jp}

\author{S. Zwicknagl}
\address{SZ:  Department of Mathematics, University of California at Riverside, 92521.}
\email{zwick@math.ucr.edu}



\begin{abstract} 
The fundamental combinatorial structure of a net in $\CP^2$ is its associated set of mutually orthogonal latin squares. We define equivalence classes of sets of orthogonal Latin squares by label equivalences of the lines of the corresponding net in $\CP^2$. Then we count these equivalence classes for small cases. Finally, we prove that the realization spaces of these classes in $\CP^2$ are empty to show some non-existence results for 4-nets in $\CP^2$.
\end{abstract}

\maketitle


\section{Introduction} \label{section:one}

In this paper we examine pairs of orthogonal Latin squares and use them to derive some non-existence results for $(4,k)$-nets in $\CP^2$. Informally speaking, an $(n,k)$-net in $\CP^2$ consists of $n$ different sets of $k$ lines in $\CP^2$ such that for every intersection point, $p$, of lines from different sets there must be exactly one line from every set passing through $p$. More precisely, we will use the following definition from \cite{Yuz}. 

 \begin{definition}\label{netdef}
A \emph{$(n,k)$-net in $\CP^2$} for $n\geq 3$ consists of a set of lines $\A \subset \CP^2$ and a finite set of points $\chi \subset \A$ such that $\A$ can be partitioned into $n$ subsets $\A =\bigcup\limits_{i=1}^n\A_i$ where $|\A_i|=k$ for all $i=1,\ldots,n$ subject to the following conditions

\begin{enumerate}
\item If $\ell_1\in\A_i$ and $\ell_2\in\A_j$ then $\ell_1\cap \ell_2\in \chi$ whenever $i\neq j$. 
\item For every $X\in \chi$ and every $i\in \{1,\ldots ,n\}$ there is exactly one line $\ell \in \A_i$ such that $X\in \ell$.
\end{enumerate}

\end{definition}

 Nets have appeared in several different areas of mathematics over the last century.  Reidemeister was one of the first to examine nets in his research on webs and their relationship to groups (see \cite{Rei}). The existence of nets has also been shown to influence the existence of finite projective planes (see \cite{DK}). More recently, Libgober and Yuzvinsky
used nets in their investigation of local systems on the complement of a complex hyperplane arrangement in $\CP^2$  (see \cite{LY}). In \cite{LY}, Libgober and Yuzvinsky investigated these nets and showed that  $(n,k)$-nets can only exist in $\CP^2$ if $n\le 5$. Then, in \cite{Yuz}, Yuzvinsky classified certain classes of nets in $\CP^2$.  In Problem 1 of \cite{Yuz} Yuzvinsky asks if there are any 4-nets in $\CP^2$ other than the Hessian configuration.  In this paper we give a partial answer to this question.  After the writing of this paper the authors learned of the overlapping but independent work of Stipins \cite{Stip, Stip1} and Urzua \cite{Urz, Urz1}.  Stipins main result in \cite{Stip} provides a more complete result than ours and proves that there are no $(n,k)$-nets in $\CP^2$ for $n \geq 4$ and $k \geq 4$.  Stipins uses geometric methods and hence provides a geometric explanation for their non-existence.  Our methods are elementary and direct.  In \cite{Urz} and \cite{Urz1} Urzua has furthered the classification of $(3,k)$-nets in $\CP^2$ by describing all possible realizations of $(3,6)$-nets and their associated moduli spaces.

Our approach towards this question uses the well known fact that one can associate $n-2$ mutually orthogonal $k\times k$ Latin squares to a $(n,k)$-net (see e.g. \cite{DK}). Recall that a Latin square is a $k\times k$
array containing the numbers $1, \ldots, k$ so that there are no repetitions of
any number within the same row or column. Now, given a $(n,k)$-net we label the $n$ sets of lines by the numbers $1,2,\ldots ,n$ and label the lines in each of the sets  by the numbers $1,2,\ldots, k$. With this labeling define $k\times k$-arrays $L^{m}$ (for $m = 1, \ldots, n-2$) by setting the $i^{th}$ entry in the $j^{th}$ column of $L^{m}$ to be $L_{ij}^{m} =\ell$, if the $i^{th}$ line of the first set and the $j^{th}$ line of the second set meet the $\ell$ line of the $m^{th}$ set in their intersection point. Every intersection of a line in the first set with a line in the second set is contained by exactly one line from the $m^{\text{th}}$ set. Thus, $L^m$ is a Latin square.   
Recall also that two $k\times k$ Latin squares  $L^1$ and $L^2$ are orthogonal if for each pair $(\ell,m)\in \{1,\ldots k\}^2$ there exists exactly one pair $(i,j)\in \{1,\ldots k\}^2$ such that $L^1_{ij}=\ell$ and $L^2_{ij}=m$. When $n\geq 4$, every intersection of a line in the first set with a line in the second set is contained in exactly one line from the $m^{\text{th}}$ set, and one line from the $m'^{\text{th}}$ set. Thus, the Latin squares $L^m$ and $L^{m'}$ are orthogonal. However, this labeling of the lines is not unique; permuting the labels of the lines in one of the sets, or the numbering of the sets can lead to different sets of orthogonal Latin squares. 

In this paper we propose the following program to classify $(n,k)$-nets:

\noindent (1) Define an equivalence relation on the set of $(n-2)$-tuples of mutually orthogonal $k\times k$-Latin squares identifying those tuples which are obtained from the same net.

\noindent (2) Choose a representative for each equivalence class and investigate whether it can be realized as a net in $\CP^2$. One obtains a system of equations whose solutions describe the moduli space of isomorphism classes of nets in $\CP^2$.     

In the present paper we will apply this approach to the case $n=4$. Denote the set of all possible pairs of orthogonal $k\times k$ Latin squares by $OLS_k$. We define the equivalence relation $\sim'$ in the next section.  We will call the set of these equivalence classes corresponding to Step (1) of our program  $OLS_k/\sim'$ . We obtain the following result.

\begin{theorem}\label{theorem:main} (see Theorems \ref{theorem:exist3}, \ref{theorem:exist4}, \ref{theorem:exist5}) Using the notation established above, we have
\begin{enumerate}
\item $|OLS_3/\sim' | = 1$.
\item $|OLS_4/\sim' | = 1$.
\item $1\le|OLS_5/\sim' | \le 2$.
\end{enumerate}
\end{theorem}

Then we calculate the realization space, as a net in $\CP^2$, of a representative of each equivalence class in Theorem \ref{theorem:main} to conclude the following theorem.

\begin{theorem}
\label{th:realizations}  The following is a complete classification of $(4,k)$-nets in $\CP^2$ up to projective isomorphism, for $k = 3, 4, 5, 6.$
\begin{enumerate}
 \item The Hessian arrangement is the only $(4,3)$-net in $\CP^2$ up to projective isomorphism.
 \item  There are no $(4,k)$-nets in $\CP^2$ for $k=4,5,6$.
\end{enumerate}
\end{theorem}
The  solution to Euler's well known  ``36 Officer Problem''  shows $|OLS_6| = 0$ (see \cite{Ta}), and so the case $k=6$ in Assertion 2 of Theorem \ref{th:realizations} is obvious. 

\begin{remark}
Many other authors have studied many different classes of orthogonal Latin squares (see \cite{DK}, \cite{CRC}, or \cite{Moor} for example). Though the literature is vast, the equivalence classes described herein appear to have not been studied previously.  Some of our results overlap with previously known results, but we include our own proofs in this new context to provide a more complete and self-contained treatment of this subject.
\end{remark}


\section{Preliminaries}
\label{se:prelim}
We will use the following notation throughout the paper.  Let $\sym_k$ be the
standard permutation group of $k$ objects.  Let $\mathcal{L}_k$ be the set of all
Latin squares of size $k$.  An element $L=\{L_{ij}\}_{i,j=1}^k \in \mathcal{L}_k$ is an $k\times k$
array containing the numbers $1, \ldots, k$ so that there are no repetitions of
any number within the same row or column.  Thus, we can associate to each Latin square and each pair $i\ne j$  the 
permutation
$\sigma_{ij}^L \in
\sym_k$ defined by
$\sigma_{ij}^L(L_{ip}) = L_{jp}$ for $p = 1, \ldots, k$.  That means $\sigma^L_{ij}$  is the permutation that sends the $p^{th}$ entry of the $i^{th}$ row to the $p^{th}$ entry of the $j^{th}$ row; we will supress the superscript  except in cases where there could be ambiguity (see Proposition \ref{theorem:perm}). A latin square is therefore uniquely described by its first column and the permutations $\sigma_{1,2}, \sigma_{2,3},\ldots, \sigma_{k-1,k}$.  For this paper, we may refer to
a Latin square by these permutations as $L(\sigma_{1{\blue ,}2}, \ldots,
\sigma_{k-1,k})$.  Since there are no repetitions in a row, the
$\sigma_{ij}$ must be fixed point free.

There are several operations of Latin
squares that preserve the property of being Latin, i.e., they define bijective maps from $\mathcal{L}_k$ to itself. Let $L \in
\mathcal{L}_k$.

\begin{tabular}{r l}
(S1):& Exchange row $i$ with row $j$. \\
(S2):& Exchange column $i$ with column $j$. \\
(S3):& Exchange two of the symbols in $L$. \\
\end{tabular}

We define a relation $\sim$ on $\mathcal{L}_k$ as follows, $L\sim L'$ if and
only if one can change $L$ into $L'$  by a finite sequence of
applications of (S1), (S2), and (S3).  

The following set, which is just the set of Latin squares that are multiplication tables of cyclic groups, will be of principle use in this note:
$$\begin{array}{c c l}
\mathcal{G}'_k & := & \{L \in \mathcal{L}_k\, | \sigma_{i,i+1}^L = \sigma_{j, j+1}^L
\text{ for all } i,j,  \text{ and } \sigma_{i,i+1}^L \text{ is a }
k-\text{cycle}\}, \\
\mathcal{G}_k & := & \{L \in \mathcal{L}_k\, | \exists H \in \mathcal{G}'_k \text{ so that } L \sim H \}\,.
\end{array}
$$
For convenience, we write an element $L \in \mathcal{G}'_k$ as $L_\sigma$, where
$\sigma$ is the associated permutation. 

The Latin squares $L$ and $L'$ are said to be
\emph{orthogonal} if the map $(i,j)\mapsto (L_{ij}, L'_{ij})$ is surjective (and
hence bijective) as a function from $\{1, \ldots, k\}\times \{1, \ldots, k\}$ to
itself.  Denote the set of all orthogonal Latin squares of size
$k$ as
$OLS_k$.  The following useful theorem is well known and gives an equivalent condition to a pair of Latin squares being orthogonal (see \cite{DK}).  A \emph{transveral} on a $k\times k$ Latin
square is a collection of $k$ entries subject to two conditions:  no
two entries are in the same row or column, and there is no repetition of the values of
the entries.

\begin{theorem}[\cite{DK}]\label{transversals}
If $L \in \mathcal{L}_k$, then there exists an $L' \in \mathcal{L}_k$ so that
$(L,L') \in OLS_k$ if and only if
there exist $k$ disjoint transversals on $L$.
\end{theorem}
Similar to the relations on Latin squares, there are several operations one can preform on pairs of
orthogonal Latin squares which preserve the property of being orthogonal.  We
list some of them below.  Suppose $(L, L') \in OLS_k$.

\begin{tabular}{r l}
(R1): & Exchange row $i$ with row $j$ in both $L$ and $L'$. \\
(R2): & Exchange column $i$ with column $j$ in both $L$ and $L'$. \\
(R3): & Exchange two symbols of $\{1, \ldots, k\}$ in $L$. \\
(R4): & Exchange two symbols of $\{1, \ldots, k\}$  in $L'$. \\
 (R5):& Transpose either $L$ or $L'$. \\
(R6): & Apply the map $(L, L') \mapsto (L', L)$.
\end{tabular}

We define a relation $\sim'$ on $OLS_k$ as follows, $(L_1,
L_2) \sim' (L_1', L_2')$ if and
only if one can change $(L_1, L_2)$ into $(L_1',
L_2')$ by a finite sequence of applications of (R1)--(R6).  The next corollary is an obvious consequence of Theorem \ref{transversals}.

\begin{corollary}  \label{theorem:cor}
If $(L, H) \in OLS_k$ and $L\sim L'$, then there exists an $H' \in \mathcal{L}$
so that $(L,H) \sim' (L', H')$.
\end{corollary}

Most authors have attempted to determine $|OLS_k|$
for various $k$ using design theory; see, for example,
\cite{Ab, Erd} and \cite{Li}.  In addition to computing their size the authors wish to understand these sets more
explicitly.  More specifically, we wish to
compute the size of the set $OLS_k/\sim'$ and to exhibit representatives for
these sets.


\section{Equivalence Classes of Orthogonal Latin Squares} \label{section:two}

In this section we explicitly describe the structure of the sets  $OLS_k/\sim'$
for $k = 3, 4$ and $5$ in Theorems \ref{theorem:exist3}, \ref{theorem:exist4}, \ref{theorem:exist5}.  In order to prove these results we first have to establish some elementary results about the sets $\mathcal{L}_k/\sim$ and $\mathcal{G}_k/\sim$.

\begin{lemma}  \label{lemma:2-1}
Adopt the notation given in Section \ref{section:one}. Fix an integer $k \geq
3$.  Let $L \in \mathcal{G}_k$, and let $L \sim L_\sigma \in \mathcal{G}_k'$.
\begin{enumerate}
\item  $L_\sigma \sim L_{(1 \cdots k)}$.  Thus, $\mathcal{G}_k/\sim$ contains
one element.
\item  Let $L_\sigma \in \mathcal{G}_k$.  There exists
one transversal to
$L_\sigma$ if and only if there exists $k$ disjoint transversals to $L_\sigma$. 
\item  If $k$ is odd, and $L \in \mathcal{G}_k$ is given, then there exists an
$L' \in \mathcal{G}_k$ which is orthogonal to $L$.
\item  If $k$ is even, and $L \in \mathcal{G}_k$, then there does not exist an
orthogonal mate to $L$.
\end{enumerate}
\end{lemma}

\begin{proof} Let $L_{\sigma} \in \mathcal{G}_k/\sim$ be given. Since $\sigma$ is a
cycle, it acts transitively on the set $\{1, \ldots, k\}$.  Thus by a
rearrangement of rows, we may order the first column from top to bottom as
$\sigma^1(1), \sigma^2(1), \ldots, \sigma^k(1)$.  Since the second
column is produced by application of $\sigma$ to each element of the first
column, the second column must now read $ \sigma^2(1), \ldots, \sigma^k(1),
\sigma^{k+1}(1) = \sigma^1(1)$.  Now relabel the symbols according to the rule
$\sigma^{p}(1) \mapsto p$.  Assertion 1 is now established.  

We now prove assertion 2.  
Let  $L_{\sigma} \in \mathcal{G}_{k}/\sim$.  Since permuting rows preserves $\sigma$,  we may assume without loss of generality that the transveral
is the main diagonal $\{(L_\sigma)_{ii}\}_{i=1}^k$.  Reading off these elements
from the top left of the square, the transversal must be, for some $j$, the
ordered set of numbers $\sigma^{j+1}(1), \sigma^{j+2}(1), \ldots,
\sigma^{j+k}(1)$.    By applying $\sigma$, we shift each location in the
transversal to the right (the element $(L_\sigma)_{kk} \mapsto (L_\sigma)_{1k}$).  This
is a new transversal, and the process may be repeated to yield $k$ disjoint
transversals.  This proves Assertion 2.

We now produce a transversal on $L_{(1\cdots k)}$ when $k$ is odd to prove
Assertion 3.  Corollary \ref{theorem:cor} allows us 
to consider only this case.  It is easy to verify that the main diagonal of
$L_{(1\cdots k)}$ is a transversal with entries $1, 3, \ldots, k, 2, 4, \ldots,
k-1$.  The hypothesis that $k$ is odd is necessary here. 
If $k$ were even the sequence down the main diagonal repeats the odd numbers
less than $k$.    Assertion 3 is now complete.

We finish the proof of Lemma \ref{lemma:2-1}
by establishing Assertion 4.  
 By Corollary \ref{theorem:cor}, it suffices to show that there does not exist a single transversal on
$L_{(1\cdots k)}$.  By Assertion 1 we may assume $\sigma:= (1\cdots k)$.    Suppose to
the contrary that there exists a transversal.  Rearrange the columns so that the
transversal is the main diagonal, noting that $\sigma_{1j} = \sigma^{i_j}$ where, $\{i_j\}$ is some ordering of the
numbers $0,
\ldots, k-1$.  Recall that $\sigma_{1j}$ is the permutation sending column 1 to column
$j$.  Notice that $i_1 = 0$ and that the numbers down the new first column
are still numbered $p, \sigma(p), \ldots,
\sigma^{k-1}(p)$ for some $p \in \{1, \ldots, k\}$.  Thus, the
numbers in the transversal are as follows:
$$
p, \sigma^{1+i_2}(p), \sigma^{2+i_3}(p), \ldots, \sigma^{k-1+i_k}(p)\,.
$$

But $\sigma^{1+i_2}(p) = \sigma^{i_2}(\sigma^1(p)) =
\sigma^{i_2}([p+1]) = [p+1+i_2]$, where $[ \cdot ]$ denotes reduction mod $k$ (reducing mod $k$ changes the numerical entries from $\{1,
\ldots, k\}$ to  the numbers $\{0, \ldots, k-1\}$, but this is irrelevant to our proof). 
Similar to above, we have $\sigma^{j-1+i_j}(p) = [p+j-1+i_j]$. Since this is a transversal,
the collection of these numbers (mod $k$) must be in set bijection with
$\mathbb{Z}_k$.  Thus the sum of all of these elements must be congruent
to $k(k-1)/2 \mod k$.  

We compute:
$$
\sum_{j = 1}^k \left[
p+j-1+i_j \right]
=
\left[\sum_{j = 0}^{k-1} j\right]  
=
 \left[ \sum_{j = 1}^{k-1} j\right] 
=
\left[ \frac{k(k-1)}{2}\right]\,.
$$
Breaking up the left hand side (notice the change of index in the middle sum),
we conclude:
$$
[kp]+ \left[\sum_{j=1}^{k-1} j\right] + \left[\sum_{j=1}^k i_j\right] =
\left[\frac{k(k-1)}{2}\right]\,.
$$
The middle sum on the left is equal to the right hand side, and $kp =0 \mod k$.  So we conclude
\begin{equation}  \label{equation:hammer}
\sum_{j=1}^k i_j = 0 \mod k.
\end{equation}

Since collectively the $i_j$ simply reorder the
numbers 1 through k,
their sum must be:
\begin{equation}  \label{equation:hammer2}
\left[\sum_{j=1}^k i_j\right]
 = \left[\sum_{j=1}^{k} j \right] =
\left[\sum_{j=0}^{k-1} j\right] = \left[\frac{k(k-1)}{2}\right]\,.
\end{equation}
But if $k$ is even, then $\frac{k(k-1)}{2}$ is not congruent to zero
mod $k$.  By Equations (\ref{equation:hammer}) and (\ref{equation:hammer2}), we have a contradiction, and it completes the proof of Lemma \ref{lemma:2-1}.
\end{proof}

\begin{theorem} \label{theorem:golden} Let $k\ge 3$.
\begin{enumerate}
\item The set $\mathcal{G}_k/\sim$ contains exactly one element.
\item The element $\mathcal{G}_k/\sim$ has an orthogonal mate if and only if $k$ is odd. 
 \end{enumerate}
\end{theorem}

\begin{proof}
 Theorem \ref{theorem:golden}  (1)  follows from Lemma \ref{lemma:2-1}  (1), since $L_\sigma\sim L_{(1\ldots k)}$ immediately yields that $\mathcal{G}_k/\sim$ contains exactly one element as asserted. Lemma \ref{lemma:2-1}  (3) and (4) establish Theorem \ref{theorem:golden}(2).
  \end{proof}

There is a stronger version of Assertion 1 of Lemma \ref{lemma:2-1} that will be of use.  Recall
that we can associate to each permutation $\sigma \in \sym_k$, a multi-index $I(\sigma)$
that describes the size and number of disjoint cycles in $\sigma$ when
it is written uniquely as a product of disjoint cycles (excluding the fixed
points of $\sigma$).   For instance, the permutation
$(15)(236)(49)
\in S_9$ corresponds to the multi-index $I((15)(236)(49))=(2,2,3)$ since there are two disjoint
$2$-cycles, and one $3$-cycle.  The following result shows
that this multi-index is the only relevant information when considering permutations in
equivalence classes of Latin squares.

\begin{proposition}  \label{theorem:perm}
Let $L \in \mathcal{L}_k$.  For each pair $(i,j)$, $i,j\in 1,\ldots,k$ there exists $W\in \mathcal{L}_k$ such that $L\sim W$ and $I(\sigma_{ij}^L)=I(\sigma_{12}^W)$.
 Moreover, one can choose $W$ such that $W_{1j} = j$.
\end{proposition}

\begin{proof}
We begin by exchanging columns $i \leftrightarrow 1$, and $j \leftrightarrow
2$.  It is easy to verify that if $\tau\in \sym_k$ is a global relabeling of the
entries of $L$ and we produce the equivalent square $W$ from such a relabeling, then
$\sigma^W_{12} = \tau \sigma^L_{ij} \tau^{-1}$.  Since the conjugation action of
$S_k$ on itself is a transitive action among permutations of the same type, the
first part of the Lemma follows.   Exchanging rows does not change the
permutation $\sigma_{12}$, thus we may reorder rows so that $W_{1j} = j$.
\end{proof}

We can now lay the framework for our study of
$4$-nets in $\CP^2$ below by separately establishing Assertions 1, 2, and 3 of Theorem \ref{theorem:main}.

\begin{theorem}  \label{theorem:exist3}
The set $OLS_3/\sim'$ contains only one element: $ (L_{(123)}, L_{(132)})$.
\end{theorem}
\begin{proof}
It is a basic fact \cite{Br} that the maximum number of mutually orthogonal Latin squares of order 3 is 2.  The only Latin squares of size
3 must belong to $\mathcal{G}_3$.  The only 3-cycles to generate these squares
are $(123)$ and $(132)$.  
\end{proof}

For the next theorem, set $\tau_1 = (12)(34), \tau_2 = (14)(23),$ and $\tau_3 =
(13)(24)$  as elements of $\sym_4$, and set $L_1:=L(\tau_1, \tau_2, \tau_1)$,
$L_2:= L(\tau_2, \tau_3, \tau_2)$, and $L_3:=L(\tau_3, \tau_1, \tau_3)$ as
$4\times 4$ Latin squares.  

\begin{theorem} \label{theorem:exist4}
Let $L_i$ and $\tau_i$ be as above.
 \begin{enumerate} 
 \item  $L_1, L_2,$ and $L_3$ are mutally orthogonal.
 \item  Any pair $(L_i,L_j)\in OLS_4/\sim'$ with $i \neq j$ are equivalent.
\item  $OLS_4/\sim' \ =  \{(L_1, L_2)\}$ contains only one element.
\end{enumerate}
\end{theorem}

\begin{proof}
One can easily verify Assertion 1.  Although the relation of orthogonality of  
Latin squares is not reflexive or transitive, it is symmetric.  Thus, to 
establish Assertion 2, we  produce a sequence of steps showing $(L_1,
L_2)
\sim' (L_2, L_3) \sim' (L_3, L_1)$.  To show each relation, simply use (R2)
and cycle (column 2 $\to$ column 3 $\to$ column 4
$\to$ column 2) in both squares.

We now prove Assertion 3.  Let $(L, L') \in OLS_4$, as in the proof of
Theorem \ref{theorem:exist3}, express $L$ in terms of its associated
fixed point free permutations $\sigma_{12}$, $\sigma_{23}$, $\sigma_{34}$, and
$\sigma_{41}$.  The fixed point free permutations in $\sym_4$ are exactly
$\tau_1, \tau_2$ and $\tau_3$, and the $4$-cycles.    We consider two cases: 
either any of $\sigma_{12}, \sigma_{23},$ or $\sigma_{34}$ are
a
$4$-cycle $\sigma$, or none of them are.

By Proposition \ref{theorem:perm}, we may assume the entries
$L_{1j} = j$, and $\sigma = (1234)$.  Now either $L_{13} = 3$ or $4$, but after
filling in the rest of the square, either choice shows $L \sim L_{(1234)}$.  By
Theorem \ref{theorem:golden} and Corollary \ref{theorem:cor}, we know that $L$
must not have an orthogonal mate, contradicting our assumption that $(L,L') \in
OLS_4$.  We conclude that each of $\sigma_{12}, \sigma_{23},$ and $\sigma_{34}$
are the permutations $\tau_1, \tau_2,$ and $\tau_3$.  

By considering the distinct permutations $\sigma_{12}, \sigma_{13},$ and
$\sigma_{14}$ instead (which must collectively be, in some order, the
permutations $\tau_1, \tau_2$ and $\tau_3$), we see that up to a change of
columns that $L \sim L_1$.  It is an easy exercise to
show that there are two possibilities for orthogonal mates: $(L_1, L_2)$ and
$(L_1,L_3)\in OLS_4$.  Assertion 3 now follows by Assertion 2.
\end{proof}

We conclude our study of orthgonal Latin squares with the following result. 

\begin{theorem}  \label{theorem:exist5}
Let $L \in \mathcal{L}_5$. 
\begin{enumerate}
\item  If $\sigma^L_{12} = (12)(345)$, then there do not exist $5$ disjoint
transversals.
\item  If $L$ has an orthogonal mate, then $L \sim L_{(12345)}$. 
\item  The possible orthgonal mates to $L_{(12345)}$ are $L_{(15432)},
L_{(14253)}$, and $L_{(13524)}$.
\item  The set $OLS_5/\sim'$  contains at most two elements: $ (L_{(12345)}, L_{(15432)})$ and $(L_{(12345)},
L_{(14253)}).$
\end{enumerate}
\end{theorem}

\begin{proof}  

Assertion 1 follows from a routine check for transversals.  One can show that if $\sigma_{12} = (12)(345)$, then there are not 5 disjoint transversals in $L$.  In fact, with the correct choices one only needs consider two disjoint transversals.

By Proposition \ref{theorem:perm}, we need to only consider two permutations: 
$\sigma_{12} = (12)(345)$ or $(12345)$.  These are the only
possible types of fixed point free permutations in $\sym_5$.  Having proved Assertion 1 if $L$ has an orthogonal mate then we may assume that $\sigma^L_{12} = (12345)$. Then Assertion 2 follows from using the assumption that $L$ is Latin and has 5 disjoint transversals. There are three different Latin squares equivalent to $L_{(12345)}$:  $L_{(15432)}, L_{(14253)}$, and $L_{(13524)}$. Further, any orthgonal mate is
also determined as one of the above squares.  Thus
Assertions 2 and 3 follow.

To prove the last assertion, we simply note that by permuting columns in the pair
$(L_{(12345)}, L_{(14253)})$ and relabeling each square seperately, we see that
$$(L_{(12345)}, L_{(14253)}) \sim' (L_{(12345)}, L_{(13524)})\,.$$  This completes
the proof of Assertion 4.
\end{proof}

\begin{remark}
Let $L = L(\sigma_{12}, \ldots, \sigma_{k-1,k}) \in \mathcal{L}_k$.  We have
proved in this section that if $k = 3,4,$ or $5$, and if $L$ is to have an
orthogonal mate, then each of the associated permutations to $L$ must be even
(i.e., expressible as a product of an even number of transpositions).  However, for higher $k$ there are examples where one of the associated permutations is odd and it has an orthogonal mate (see \cite{DK, Moor}). There are also examples where all the permutations are even, but there is no orthogonal mate.  It would be interesting to know how the parity of the $\sigma_{ij}$ effects the existence of an orthogonal mate.
\end{remark}


\section{Realization Spaces} \label{section:three}

As discussed in the introduction, there is a relationship between $(4,k)$-nets
and $OLS_k/\sim'$.  Given a
$(4,k)$-net, one can construct a unique element of $OLS_k/\sim'$ that
represents the underlying combinatorial structure.  Nets in $\CP^2$ that are projectively isomorphic will produce equivalent combinatorial structures.  Conversely, given an
element of $OLS_k/\sim'$ there need not exist a
$(4,k)$-net in $\CP^2$ with the given structure.

We quickly review how to compute the realization space in $\CP^2$ of a pair of orthogonal Latin squares (see \cite{BLSWZ}). A pair of Orthogonal Latin squares $(L_1,L_2)$ defines the points of $\chi$ for a combinatorial structure of a $(4,k)$-net that might be associated with a $(4,k)$-net in $\CP^2$. Let $M_{(L_1,L_2)}$ be a $4k \times 3$ matrix of complex numbers, defined by $4$ blocks of $k$ rows where the $i^{th}$ row in the $j^{th}$ block is $a_{j,i} \ b_{j,i} \ c_{j,i}$. Let the rows of $M_{(L_1,L_2)}$ be the coefficients of the linear forms defining the lines of the alleged $(4,k)$-net in $\CP^2$. Then for each point of $\chi$, the corresponding minor of $M(L_1,L_2)$ should be zero. Thus, the realization space of the pair of orthogonal Latin squares $(L_1,L_2)$ is the space of solutions to all of the minors associated to $\chi$; we denote this space $R(L_1,L_2)$. Each line of the net to be realized is labeled by a distinct element of $\{1,\ldots ,4k\}$; hence the points of $\chi$ are given by 4-tuples of distinct elements of this set. In each case we use the lexicographic ordering of the 4-tuples to compute the minors consecutively.
The next proposition shows that we only need to realize one representative of each equivalence class. This proposition is a consequence of the following two facts: 1) the relations (R1)-(R6) preserve the isomorphism type of the Latin squares component of the intersection lattice of a possible net in $\CP^2$ and 2) the realization space of a lattice or matroid is invariant under isomorphism.

\begin{proposition}
\label{theorem:netequivalences}
 If $(L_1,L_2) \sim'(L_1',L_2')$ then $R(L_1,L_2)$ and $R(L_1',L_2')$ are isomorphic as varieties.
\end{proposition} 

\begin{remark} Using Proposition 3.3 and Corollary 3.5 in \cite{Yuz} we can assume that \begin{equation*}B= \left[\begin{array}{ccc}a_{1,1}&b_{1,1}&c_{1,1}\\
a_{1,2}&b_{1,2}&c_{1,2}\\
a_{1,3}&b_{1,3}&c_{1,3}\\
a_{2,1}&b_{2,1}&c_{2,1}\end{array}\right]=\left[\begin{array}{ccc}1&0&0\\ 0&1&0\\ 0&0&1\\ 1&1&1\end{array}\right]_.\end{equation*}
\end{remark}

\subsection{Realization of (4,3)-nets/Proof of Theorem 1.3 (1)}\label{realize3}

Theorem
\ref{theorem:exist3} shows that the only combinatorial structure possible for a $(4,3)$-net is given by the pair  $(L_{(123)},L_{(132)})\in OLS_3/\sim'$.  By computing minors, we conclude
$$M_{(L_{(123)},L_{(132)})}=\left[\begin{array}{ccc}
&B&\\
1&\w^i&\w^j
\end{array}\right] $$ 
where $\w$ is a (primitive) root of $x^3-1$ and $i$ and $j$ range through the set $\{ 0,1,2\}$. This is the Hessian configuration, see Example 6.29 of \cite{OT} and Example 3.6 of \cite{Yuz}.   This proves Theorem \ref{th:realizations} for $k=3$. \hfill $\qedbox$


\subsection{Realization of (4,4)-nets/Proof of Theorem 1.3 (2)}\label{realize4}

By Theorem \ref{theorem:exist4} and Proposition \ref{theorem:netequivalences}, there is only one combinatorial
structure for a $(4,4)$-net.  We prove Theorem \ref{th:realizations} in the
case $k = 4$ by attempting to realize the combinatorial structure given by $(L_1,
L_2) \in OLS_4/\sim'$, where $L_1$ and $L_2$ are the squares given in Theorem
\ref{theorem:exist4}. In this case $|\chi |=16$. Using the first $13$ points of $\chi$ given by the pair $L_1$ and $L_2$ we get that $$M_{(L_1,L_2)}=\left[ \begin{array}{ccc}
&B&\\
r-t^{-1}&1+r&1-r\\
t&-t&1\\
t&1&-1\\
-t^{-1}&1&t^{-1}\\
t&1&1\\
1&-t&1\\
1&1&-1\\
-1&1&t^{-1}\\
-t^{-1}&1&1\\
-1&-t&1\\
1&-1&1\\
1&1&t^{-1}
\end{array}\right]_.$$ Then using the last three points of $\chi$ we find that $r=\frac{1}{4}(1-t^{-2})$, $t=-2\pm \sqrt{5}$, and $t^2+3=0$.  This system of equations has no solution, so $R(L_1,L_2)=\emptyset$ proving Theorem \ref{th:realizations} for $k = 4$: there do not exist any $(4,4)$-nets in $\CP^2$.  \hfill $\qedbox$


\subsection{Realization of (4,5)-nets/Proof of Theorem 1.3 (3)} 

By Theorem \ref{theorem:exist5} and Proposition \ref{theorem:netequivalences}, there are at most two possible combinatorial
structures for a $(4,5)$-net. In this case $|\chi |=25$. First, we consider the pair $(L_{(12345)}, L_{(15432)})$ as
our combinatorial structure. Using only the first $17$ points of $\chi$, we find that for some number $t$, there is a line in class 2
corresponding to the line defined by the linear form $tx + y + t^2 z = 0$, and a line in class 3 that corresponds to the line defined by the linear form $x + y + tz = 0$.  Then using one more point of $\chi$ we get that $t = 0$ or
$1$.  This cannot happen, in both cases a line repeats. Hence the realization space is empty and there does not exist any $(4,5)$-net in $\CP^2$ with Latin squares $(L_{(12345)}, L_{(15432)})$.

Now, we compute the realization space for the pair of Latin squares $(L_{(12345)},L_{(14253)})$. Using the first $22$ points of $\chi$ we compute that $$M_{(L_{(12345)},L_{(14253)})}=\left[\begin{array}{ccc}
&B&\\
1-qt^{-3}&1-q&1-qt^{-1}\\
a&b&c\\
t^{-1}&t^{-3}&1\\
t^{-1}&1&t\\
t&1&t^{-1}\\
t&t^3&1\\
t&1&1\\
1&t^{-3}&1\\
1&1&t\\
t^{-3}&1&t^{-1}\\
t^2&t^3&1\\
t^{-1}&1&1\\
t^{-2}&t^{-3}&1\\
t^3&1&t\\
1&1&t^{-1}\\
1&t^3&1
\end{array}\right]$$ where $t^5-1=0$ and $q=\frac{1+t^3}{2+t^3-t^2}.$ Then using one more point of $\chi$ we get that $t^3-1=0$. Hence, $R (L_{(12345)},L_{(14253)})$ and there does not exist any $(4,5)$-nets in $\CP^2$.  This concludes the proof of Theorem \ref{th:realizations}.  \hfill $\qedbox$

The second author and Eric Dybeck have made further progress towards writing a computer program that would determine the classification of $(4,k)$-nets in $\CP^2$ for $k \geq 7$.

\section{Acknowledgements}
The first author is partially supported by a CSUSB faculty grant. The third author has been supported by NSF grant \# 0600893 and the NSF Japan program.
The authors would like to thank Sergey Yuzvinsky for many helpful conversations.


\end{document}